\documentclass[12pt,twoside,a4paper]{amsart}
\usepackage{amsmath,amssymb}
\usepackage{amsfonts}
\usepackage{mathrsfs} 
\usepackage{mathtools}
\usepackage{cite}
\usepackage{graphicx}
\usepackage{txfonts}
\usepackage{mathrsfs}
\usepackage{empheq}
\newcommand{\ds}{\displaystyle}

\newcommand{\an}{E^{\text{AN}}(x)}
\newcommand{\AS}{E_{\sigma_1}^{\text{AR}}(x)}
\newcommand{\as}{E_{\sigma_1}^{\text{AN}}(x)}
\newcommand{\e}{\varphi(n)}
\newcommand{\si}{\sigma_1 (n)}
\newcommand{\s}[2]{\{ #1 (#2) \}_{#2 = 1}^\infty}

\allowdisplaybreaks[4]
\theoremstyle{definition}
\newtheorem{definition}{Definition}[subsection]
\newtheorem{remark}[definition]{Remark}

\theoremstyle{theorem}
\newtheorem{theorem}{Theorem}[subsection]
\newtheorem{lemma}[theorem]{Lemma}

\newtheorem{fact}[theorem]{Fact}

\numberwithin{equation}{subsection}
\mathtoolsset{showonlyrefs=true}
\mathtoolsset{showonlyrefs=false}

\begin{document}

\title[ A boundary example $\ldots$ ]{A boundary example for the Volterra decomposition of arithmetic error terms}
\author{Hideto Iwata}
\address{National Institute of Technology, Gunma College, Japan.}
\email{iwata@gunma-ct.ac.jp}
\keywords{Volterra integral equation of the second type, remainder term in an asymptotic formula, arithmetical function, divisor function, analytic part.}
\subjclass[2020]{11M26, 11N37, 45D05}
\maketitle

%\textbf{Abstract.}\hspace{0.1cm}
%Kaczorowski and Wiertelak introduced an arithmetic - analytic decomposition of the error term in the summatory Euler totient function $\varphi(n)$. In this paper, we establish an analogous decomposition for the divisor function $\sigma_1(n)$ via a generalized Volterra integral equation. We obtain an explicit expression for the analytic part and derive upper bounds parallel to those in the Euler totient case. Assuming the Riemann Hypothesis, we prove
%\[ {\as} \ll x^{\delta'} \exp\!\left( \frac{\log x}{\log\log x} \right) \]and, as a consequence,
%\[ {\as} \ll_\epsilon x^{\delta' + \epsilon} \qquad (\epsilon>0). \]We also explain why the full equivalence with the Riemann Hypothesis, known for $\varphi(n)$, does not extend to $\sigma_1(n)$.

\begin{abstract}
Let $\varphi(n)$ denote the Euler totient function. We study the analytic part associated with the summatory function of $\sigma_1(n)$ and obtain explicit bounds under the Riemann Hypothesis. In particular, we establish an upper bound of order
$x^{\delta'} \exp\!\left( \frac{\log x}{\log\log x} \right)$, where $\delta' = \max\{1/2,\delta\}$.
\end{abstract}

\section{Introduction and the statement of the main result}
\subsection{Introduction}

Let $\varphi(n)$ denote the Euler totient function and consider the error term
\begin{equation}
   E(x) = \sum_{n \le x} \varphi(n) - \frac{3}{\pi^2} x^2.   \label{eq:1.1.1}
\end{equation}In their work, Kaczorowski and Wiertelak \cite{Kac and Wie} studied a Volterra integral equation associated with $E(x)$ and introduced a decomposition of $E(x)$ into an arithmetic part and an analytic part. A remarkable feature of their result is that the analytic part $E^{\mathrm{AN}}(x)$ is closely connected with the Riemann Hypothesis; in fact, certain estimates for $E^{\mathrm{AN}}(x)$ are equivalent to the Riemann Hypothesis. 

The purpose of this paper is not merely to extend the construction of Kaczorowski and Wiertelak to another arithmetical function, but rather to investigate the boundary behavior of the Volterra framework established in ~\cite{I}. The general theory in \cite{I} requires the convergence of the series 

\[ \ds \sum_{n=1}^{\infty}\frac{a(n)}{n^2}. \] 

Although an arithmetic part and an analytic part can still be constructed, several phenomena differ essentially from those occurring in the classical case of the Euler totient function. In particular, the divisor function provides the first example within the framework of ~\cite{I} for which (i) an additional $O(x^{1/2})$-term appears in the decomposition of the error term, and (ii) the Riemann-Hypothesis-equivalent condition established by Kaczorowski and Wiertelak disappears. These phenomena indicate that the behavior of the analytic part can change drastically even within the same Volterra framework. The purpose of the present paper is therefore not only to obtain estimates for the analytic part associated with $\sigma_1 (n)$, but also to identify a concrete situation in which the qualitative behavior of the framework differs from the classical Euler totient case. In particular, we consider the divisor function
\[ \ds \sigma_1(n) = \sum_{d \mid n} d, \] and study the corresponding error term

\begin{equation}
   E_{\sigma_1}(x) = \sum_{n \le x} \sigma_1 (n) - \frac{\pi^2}{12} x^2.   \label{eq:1.1.2}
\end{equation}

Using a general result on Volterra integral equations developed in \cite{I}, we introduce a decomposition of $E_{\sigma_1}(x)$ into an arithmetic part and an analytic part, denoted by $E^{\mathrm{AR}}_{\sigma_1}(x)$ and $E^{\mathrm{AN}}_{\sigma_1}(x)$, respectively. Our first main result gives an explicit representation of the analytic part $E^{\mathrm{AN}}_{\sigma_1}(x)$ in terms of a Mellin transform involving the product $\zeta(s)\zeta(s-1)$. Based on this representation, we derive upper bounds for $E^{\mathrm{AN}}_{\sigma_1}(x)$ under the assumption of the Riemann Hypothesis. A key point of this paper is that the behavior of the analytic part in the case of $\sigma_1(n)$ is fundamentally different from that in the case of $\varphi(n)$. More precisely, while in the totient case the analytic part yields an equivalence condition for the Riemann Hypothesis, such an equivalence does not arise in the present setting. This difference is explained by the structure of the associated Dirichlet series:
\begin{gather*}
   \sum_{n=1}^{\infty} \frac{\varphi(n)}{n^s} = \frac{\zeta(s-1)}{\zeta(s)}, 
   \\
   \sum_{n=1}^{\infty} \frac{\sigma_1(n)}{n^s} = \zeta(s)\zeta(s-1),
\end{gather*}and in particular by the absence of the reciprocal factor $1 / \zeta(s)$ in the latter case. Thus, the present work provides a contrasting example to the results of Kaczorowski and Wiertelak, illustrating how the arithmetic nature of the underlying Dirichlet series affects the analytic behavior of the associated error terms. The paper is organized as follows. In Section 2, we construct the arithmetic and analytic parts for $E_{\sigma_1}(x)$. In Section 3, we derive a Mellin transform formula and its inverse representation. In Section 4, we prove the main theorems. A fundamental result in this direction is the following theorem due to
Kaczorowski and Wiertelak.

\begin{theorem}[Theorem 1.7 in~\cite{Kac and Wie}]
The following statements are equivalent.
\begin{enumerate}
\item[(1)] The Riemann Hypothesis is true.
\item[(2)] There exists a positive constant $A$ such that for $x \geq e^e$ we have
\begin{equation}  
   \an \ll x^{1/2}\exp \left( A\frac{\log x}{\log\log x} \right).   \label{eq:1.1.3}
\end{equation}
\item[(3)] For every $\epsilon>0$ and $x \geq 1$ we have
\begin{equation} 
   \an \ll_\epsilon x^{1/2+\epsilon}.
   \label{eq:1.1.4}
\end{equation}
\end{enumerate}
\end{theorem} This naturally leads to the question of how the behavior of the analytic part changes for other functions arising in the framework of ~\cite{I}. In this paper, we carry out the same type of analysis as in~\cite{Kac and Wie} for an error term in the asymptotic formula for $\si$ in place of $\e$. Then we obtain Theorems 1.2.1 and 1.2.3 in the next subsection, corresponding to Theorem 1.1.1.

\subsection{Main theorem}

The explicit form of $\as$ is given in \eqref{eq:2.1.14} in section 2. First, we obtain the following bound for $\as$ :
\begin{theorem}
Under the Riemann Hypothesis, we have
\begin{equation}
   \as \ll x^{\delta'} \exp\left( \frac{\log x}{\log\log x} \right),   \label{eq:1.2.1}
\end{equation}
for $x \ge e^{e}$.  
Here, $\delta'$ is defined by $\delta' = \max\{1/2, \delta\}$, where $0 < \delta < 1$ is arbitrary.
\end{theorem}

In what follows, we write a complex variable as $s = \sigma + it$, and we shall continue to use this notation throughout the paper. The proof of \eqref{eq:1.2.1} is based on Lemma 1.2.2 below. 

\begin{lemma}[Lemma 2.5 in~\cite{K}]
Assuming the Riemann Hypothesis, there exist positive constants $t_0$ and $A$ such that for
\[
   \frac{1}{2} < \sigma \le \frac{1}{2} + \frac{1}{\log\log |t|}
   \qquad (|t| \ge t_0)
\]
we have
\begin{equation}
   \frac{\zeta(s-1)}{\zeta(s)}
   \ll
   |t| \exp\left( 
      A \frac{\log |t|}{\log\log |t|}
      \log\frac{2}{(\sigma - \tfrac12)\log\log |t|}
   \right).
   \label{eq:1.2.2}
\end{equation}
\end{lemma}Lemma 1.2.2 follows from standard bounds for $\zeta(s)$ in the region $\sigma < 0$ together with classical estimates for $\zeta(s)^{-1}$ under the Riemann Hypothesis (see \cite{t}). Next, a bound for $\as$ obtained by Theorem 1.2.1 is as follows:

\begin{theorem}
Assuming Theorem 1.2.1, we have
\begin{equation}
   \as \ll_\epsilon x^{\delta' + \epsilon},
   \label{eq:1.2.3}
\end{equation}
where $\epsilon$ is any positive real number.
\end{theorem}

The equivalence in Theorem 1.1.1 provides a characterization of the Riemann Hypothesis in terms of $\an$. However, a corresponding equivalence for $\as$ in the asymptotic formula for $\si$ does not hold in general.

\section{Background}
\subsection{The divisor function as a boundary case}
In this section we explain why the divisor function $\sigma_1 (n)$ is distinguished among the examples arising from Theorem 2.1.1. In~\cite{I}, the author obtained a solution to a Volterra integral equation of the second type in the following setting:

\begin{theorem}[Theorem in~\cite{I}]
Let $\s{a}{n}$ be a complex-valued arithmetical function for which the series
\begin{equation}
   \sum_{n=1}^\infty \frac{a (n)}{n^2}
   \label{eq:2.1.1}
\end{equation}
is convergent with sum $2\alpha$, where $\alpha$ is an arbitrary complex number. Let  $\s{b}{n}$  be the arithmetical function defined by
\begin{equation}
   b (n) = \sum_{d\mid n} a (d)\frac{n}{d}.
   \label{eq:2.1.2}
\end{equation}
Assume that, as $x \to \infty$,
\begin{equation}
   \sum_{n \leq x} b (n) = M(x) + \text{Er}(x),
   \label{eq:2.1.3}
\end{equation}
where
\begin{gather*}
   M(x) := \alpha x^2,\\
   \text{Er}(x) := \sum_{n \leq x} b (n) - M(x).
\end{gather*}
Now we consider the following Volterra integral equation of the second type:
\begin{equation}
   F_1 (x) - \int_{0}^x F_1 (t)\frac{dt}{t} = \text{Er}(x) \quad (x \geq 0).
   \label{eq:2.1.4}
\end{equation}
Then, for every complex number $A$, the function
\begin{equation}
   F_1 (x) = (f (x) + A)x \quad (x \geq 0), 
   \label{eq:2.1.5} 
\end{equation}
is a solution of the integral equation \eqref{eq:2.1.4}, and these exhaust all solutions of \eqref{eq:2.1.4}. Here,
\begin{equation} 
   f (x) = -\sum_{n=1}^\infty \frac{a(n)}{n}\left\{ \frac{x}{n} \right\}
   \label{eq:2.1.6}
\end{equation}
for every $x \geq 0$.
\end{theorem}

Now, the function $f(x)$ is locally bounded and so the series on the right-hand side of \eqref{eq:2.1.6} is convergent (see~\cite{I}). In~\cite{Kac and Wie}, the arithmetical functions $\s{a}{n},\s{b}{n}$ in Theorem~2.1.1 are $\mu(n)$ and $\e$, respectively, and all the hypotheses in Theorem 2.1.1 are satisfied. Another choice of $\s{a}{n}, \s{b}{n}$ is $u(n)$ and $\si$, respectively, where the arithmetic function $u(n) = 1$ for all positive integers $n$. Then the function $f(x)$ in the case $a(n) = u(n)$ is
\begin{equation}
   f_2 (x) = -\sum_{n=1}^\infty \frac{1}{n}\left\{ \frac{x}{n} \right\}
   \label{eq:2.1.7}
\end{equation}
for every $x \geq 0$.To obtain the analytic part associated with $E_{\sigma_1} (x)$, we define the auxiliary function for $x \geq 0$ by
\begin{equation}
   R(x) = \text{Er}(x) - x f(x).
   \label{eq:2.1.8}
\end{equation}
Then we use the following lemma: 
\begin{lemma}[Lemma 1 in~\cite{I}]
For all positive $x$,
\begin{equation}
      R(x) = -\int_{0}^x f(t)\,dt.
      \label{eq:2.1.9}
\end{equation}
\end{lemma}

Using Lemma 2.1.2 for $f_2 (x)$, the right-hand side of \eqref{eq:2.1.9} is
\begin{align}
   R(x)   \notag
   &= -\int_{0}^x f_2 (t)\,dt 
   = -\int_{0}^x \left( -\sum_{n=1}^\infty \frac{1}{n}\left\{ \frac{t}{n} \right\} \right)dt
   = \sum_{n=1}^\infty \int_{0}^{x/n}  \{ u \}\,du   \notag\\
   &= \frac{1}{2}\sum_{n=1}^\infty \left\{ \frac{x}{n} \right\}^2 + \frac{1}{2}\sum_{n=1}^\infty \left[ \frac{x}{n} \right].
   \label{eq:2.1.10}
\end{align}
To obtain the right-hand side of \eqref{eq:2.1.10}, we use the following formula (see~\cite{Kac and Wie}, p.~2692): for $x>0$, 
\[
   \int_0^x \{ t \}\,dt = \frac{1}{2}\{ x \}^2 + \frac{1}{2}[x].
\]
Now, the second series involving the integer part of $x/n$ on the right-hand side of \eqref{eq:2.1.10} satisfies
\begin{equation} 
   \sum_{n=1}^\infty \left[ \frac{x}{n} \right]  
   = \sum_{n=1}^\infty \sum_{m \leq x/n} 1
   = \sum_{n \leq x} \sum_{\ell \mid  n} 1
   = \sum_{n \leq x} d(n)
   = x\log x + (2\gamma -1)x + O(x^{1/2}),
   \label{eq:2.1.11}
\end{equation}
as $x \to \infty$. Here, we use the following classical asymptotic formula for the divisor function $d(n) = \sum_{\ell \mid n}  1$:
\[
   \sum_{n \leq x} d(n) = x\log x + (2\gamma -1)x + O(x^{1/2}),
\]
where the constant $\gamma$ is Euler's constant. The appearance of the term $O(x^{1/2})$ is a direct consequence of the failure of M\"{o}bius cancellation in the divisor function setting. For $x \geq 1$, we define the function $g_2 (x)$ by the series
\begin{equation}
   g_2 (x) = \sum_{n=1}^\infty \left\{ \frac{x}{n} \right\}^2.
   \label{eq:2.1.12}
\end{equation}
Substituting \eqref{eq:2.1.11} and \eqref{eq:2.1.12} into \eqref{eq:2.1.10}, we obtain the following decomposition for $E_{\sigma_1}^{\text{AN}}(x)$:
\begin{equation}
   E_{\sigma_1} (x) =  \AS + \as + O(x^{1/2}),
   \label{eq:2.1.13}
\end{equation}
where
\begin{equation}
   \AS = xf_2 (x), \quad \text{and} \quad \as = \frac{1}{2}g_2 (x) + \frac{x}{2}(\log x + 2\gamma - 1)
   \label{eq:2.1.14}
\end{equation}
with $f_2 (x)$ and $g_2 (x)$ given by \eqref{eq:2.1.7} and \eqref{eq:2.1.12}, respectively.We call $\AS$ and $\as$ the \textit{arithmetic} and \textit{analytic} parts of $E_{\sigma_1} (x)$, respectively, in analogy with the decomposition of $E(x)$.

\begin{remark}
In~\cite{Kac and Wie}, the series
\begin{equation}
   \sum_{n=1}^\infty \mu(n) \left[ \frac{x}{n} \right]
   \label{eq:2.1.15}
\end{equation}
corresponds to the series involving the integer part of $x/n$ in \eqref{eq:2.1.10}. The series \eqref{eq:2.1.15} is calculated as follows:
\[
   \sum_{n=1}^\infty \mu(n) \left[ \frac{x}{n} \right]
   = \sum_{n=1}^\infty \mu(n) \sum_{m \leq x/n} 1
   = \sum_{n \leq x} \sum_{\ell \mid n} \mu(\ell)
   = 1.
\]
Hence, the series \eqref{eq:2.1.15} can be evaluated explicitly and $\an$ is represented without an error term, in contrast to \eqref{eq:2.1.11}.
\end{remark}

From the calculation in \eqref{eq:2.1.10}, we see that
\[
   \as = -\int_{0}^x  f_2 (t)\,dt \ll x^2,
\]
and so $f_2 (t)$ is integrable. Since $\as$ is represented as an indefinite integral of an integrable function, we see that $\as$ is an absolutely continuous function.

\section{Lemmas}
\subsection{Lemmas}
We prepare the Mellin transform and its inverse for $E_{\sigma_1} (x)$ in order to obtain Theorems 1.2.1 and 1.2.2. First, we prove the Mellin transform for $\as$. 

\begin{lemma}
For $\sigma>2$, we have
\begin{equation}
   \int_{1}^\infty  E_{\sigma_1}^{\text{AN}}(x)x^{-s-1}\,dx
   = \frac{\pi^2}{12} \frac{1}{s-2} + \frac{\zeta(s)\zeta(s-1)}{s(1-s)} + O(1).
   \label{eq:3.1.1}
\end{equation}
\end{lemma}

\textbf{Proof of Lemma 3.1.1.}
According to \eqref{eq:2.1.13} and \eqref{eq:2.1.14}, we have $\as = E_{\sigma_1} (x) - \AS + O(x^{1/2})$. Substituting this into the integral on the left-hand side of Lemma~3.1.1, we obtain
\begin{align}
   \int_{1}^\infty  \as x^{-s-1}\,dx
   &= \int_{1}^\infty E_{\sigma_1}(x) x^{-s-1}\,dx - \int_{1}^\infty f_2 (x)x^{-s}\,dx   \notag\\
   &\quad + O_\sigma \left( \int_{1}^\infty x^{-\sigma - 1/2}\,dx \right)   \notag\\
   &= \int_{1}^\infty E_{\sigma_1}(x) x^{-s-1}\,dx - \int_{1}^\infty f_2 (x)x^{-s}\,dx + O(1).
   \label{eq:3.1.2}
\end{align}First, we calculate the integral on the left-hand side of \eqref{eq:3.1.2}. By the definition of $E_{\sigma_1} (x)$ in \eqref{eq:1.1.2} and the condition $\sigma>2$, we have
\begin{equation}
   \int_{1}^\infty E_{\sigma_1}(x) x^{-s-1}\,dx
   = \int_{1}^\infty A(x)x^{-s-1}\,dx - \frac{\pi^2}{12} \frac{1}{s-2},
   \label{eq:3.1.3}
\end{equation}
where
\[
   A(x) = \sum_{n \leq x} \si.
\]
Using partial summation and \eqref{eq:1.1.2} again, the integral on the right-hand side of \eqref{eq:3.1.3} becomes
\begin{equation}
   \int_{1}^\infty A(x)x^{-s-1}\,dx 
   = \frac{1}{s}\sum_{n=1}^\infty \frac{\si}{n^s} = \frac{1}{s}\zeta(s)\zeta(s-1)
   \label{eq:3.1.4}
\end{equation}for $\sigma>2$. Here, we use the following Dirichlet series expansion for $\zeta(s)\zeta(s-1)$:
\begin{equation}
   \sum_{n=1}^\infty \frac{\sigma_1 (n)}{n^s} = \zeta(s)\zeta(s-1) \hspace{0.3cm}  (\sigma > 2).
   \label{eq:3.1.5}
\end{equation}
Using \eqref{eq:3.1.3} and \eqref{eq:3.1.4}, we have
\begin{equation}
   \int_{1}^\infty E_{\sigma_1}(x)x^{-s-1}\,dx
   = \frac{1}{s}\zeta(s)\zeta(s-1) - \frac{\pi^2 }{12}\frac{1}{s-2}
   \label{eq:3.1.6}
\end{equation}for $\sigma > 2$. On the other hand, by the definition of $f_2 (x)$ in \eqref{eq:2.1.7}, the integral whose integrand is $f_2 (x)x^{-s}$ on the right-hand side of \eqref{eq:3.1.2} is
\begin{align}
   \int_{1}^\infty f_2 (x)x^{-s}\,dx
   &= \int_{1}^\infty \left( -\sum_{n=1}^\infty \frac{1}{n}\left\{ \frac{x}{n} \right\} \right)x^{-s}\,dx   \notag\\
   &= -\int_{1}^\infty \left( \sum_{n=1}^\infty \frac{x}{n^2} \right)x^{-s}\,dx
      + \int_{1}^\infty \left( \sum_{n=1}^\infty \frac{1}{n}\left[ \frac{x}{n} \right] \right)x^{-s}\,dx   \notag\\
   &= -\frac{\pi^2}{6} \frac{1}{s-2}
      + \int_{1}^\infty \left( \sum_{n \leq x} \frac{\si}{n} \right)x^{-s}\,dx.
   \label{eq:3.1.7}
\end{align}
As in the derivation of \eqref{eq:3.1.4}, we have
\begin{equation*}
   \int_{1}^\infty \left( \sum_{n=1}^\infty \frac{\sigma_1 (n)}{n} \right) x^{-s}\,dx
   = \frac{1}{s-1}\sum_{n=1}^\infty \frac{\si}{n^s}
   = \frac{1}{s-1}\zeta(s)\zeta(s-1)
\end{equation*}
for $\sigma>2$. Substituting this into \eqref{eq:3.1.7}, we obtain
\begin{equation}
   \int_{1}^\infty f_2 (x)x^{-s}\,dx
   = -\frac{\pi^2}{6}\frac{1}{s-2} + \frac{1}{s-1} \zeta(s)\zeta(s-1).
   \label{eq:3.1.8}
\end{equation}
By substituting \eqref{eq:3.1.6} and \eqref{eq:3.1.8} into the right-hand side of \eqref{eq:3.1.2}, the conclusion of Lemma~3.1.1 follows. \qed

Next, we take the inverse Mellin transform in \eqref{eq:3.1.1}. 

\begin{lemma}
For any positive real number $0 < \delta < 1$, we have 
\begin{equation}
   \as = \frac{1}{2\pi i}\int_{\mathcal{L}} \zeta(s)\zeta(s-1) \frac{x^s}{s(1-s)}\,ds + O(x^{3+\delta}).
   \label{eq:3.1.9} 
\end{equation}
Here, the path of integration $\mathcal{L}$ consists of the half-line $[3+\delta - i\infty, 3+\delta - 3i]$, the semicircle $s = 3 + \delta + 3e^{i\theta}$, $\pi / 2 < \theta < 3\pi/2$, and the half-line $[ 3 + \delta + 3i, 3 + \delta + i\infty]$. 
\end{lemma}

Now we use the following fact for the inverse Mellin transform:

\begin{fact}[~\cite{T}]
Let $f(y)y^{k-1} \ (k>0)$ belong to $L(0,\infty)$, and let $f(y)$ be of bounded variation in a neighborhood of the point $y=x$. For $s = k +it$, let 
\begin{align*}
   \mathscr{F}(s) = \int_{0}^\infty f(x)x^{s-1}\,dx \hspace{0.3cm} (s=k+it).
\end{align*}
Then
\begin{align}
   \frac{1}{2}\{ f(x+0) + f(x-0) \}
   = \frac{1}{2\pi i}\lim_{T \to \infty} \int_{k-iT}^{k+iT} \mathscr{F}(s)x^{-s}\,ds.
   \label{eq:3.1.10}
\end{align}
\end{fact}

\textbf{Proof of Lemma 3.1.2.}
First, by the definition of $\as$ in \eqref{eq:2.1.14}, we see that $x^{-\sigma-1}\as \in L(1,\infty)$ for $\sigma>2$. Secondly, by the definition of $g_2 (x)$ in \eqref{eq:2.1.12} and of $\as$ in \eqref{eq:2.1.14} again, we can see that $\as$ is represented as the difference of two monotone increasing functions, so $\as$ is of bounded variation. Using Fact~3.1.3, we can take the inverse Mellin transform in \eqref{eq:3.1.1} and obtain
\begin{align}
   \as
   &=  \frac{1}{2\pi i} \int_{(3+\delta)} \left\{ \frac{\pi^2}{12} \frac{1}{s-2} + \frac{\zeta(s)\zeta(s-1)}{s(1-s)} + O(1) \right\}x^s \,ds   \notag\\
   &= \frac{\pi^2}{12}x^2 + \frac{1}{2\pi i}\int_{(3+\delta)} \zeta(s)\zeta(s-1) \frac{x^s}{s(1-s)}\,ds + O\left( x^{3+\delta} \right)
   \label{eq:3.1.11}
\end{align}
for $x>1$. Here, the path of integration is the vertical line $\text{Re} (s) = 3 + \delta$ in the whole $s$-plane. We calculate the integral in the second line of \eqref{eq:3.1.11}. Using the Laurent expansion for $\zeta(s)$ at $s=1$,
\[
   \zeta(s) = \frac{1}{s-1} + \gamma + O(|s-1|),
\]
the residue of the integrand at $s=1$ is $(\gamma - 1 + \log (2\pi x))x/2$. On the other hand, the residue at $s = 2$ is $-\pi^2 x^2 /12$. By the residue theorem we have
\begin{align}
   \frac{1}{2\pi i}\int_{(3+\delta)} \zeta(s)\zeta(s-1)\frac{x^s}{s(1-s)}\,ds
   &= \frac{1}{2\pi i}\int_{\mathcal{L}} \zeta(s)\zeta(s-1) \frac{x^s}{s(1-s)}\,ds   \notag\\
   &+ \frac{x}{2}(\gamma - 1 +\log (2\pi x)) - \frac{\pi^2}{12}x^2 ,
   \label{eq:3.1.12}
\end{align}
where the path of integration $\mathcal{L}$ is as in Lemma~3.1.2. Since the left-hand side of \eqref{eq:3.1.12} equals $\as - \pi^2 x^2 /12 + O(x^{3+\delta})$ by \eqref{eq:3.1.11}, we obtain the desired result. \qed

\section{Proof of the main theorems}
\subsection{Proof of Theorem 1.2.1}
In this section, we begin the proofs of the main theorems. First, we prove Theorem~1.2.1. Suppose that the Riemann Hypothesis holds and consider the rectangle
\[
   R = \{ s = \sigma + it \mid 2 + \delta \leq \sigma \leq 3 + \delta,\ -T \leq t \leq x^2,\ x^2 <T \}.
\]
Since the integral along $R$ is zero by Cauchy's integral theorem and the integral from $2 + \delta -iT$ to $3+ \delta -iT$ tends to zero as $T \to \infty$, we can shift the path of integration and obtain, for $x > 1$,
\begin{align}
   \int_{3 + \delta - i\infty}^{3 + \delta - ix^2} \frac{\zeta(s)\zeta(s-1)}{s(1-s)}x^s \,ds
   &= \left\{ \int_{2 + \delta - i\infty}^{2 + \delta - ix^2}
      +  \int_{2 + \delta - ix^2}^{3 + \delta - ix^2} \right\}
      \frac{\zeta(s)\zeta(s-1)}{s(1-s)}x^s \,ds   \notag\\
   &=  \int_{2 + \delta - ix^2}^{3 + \delta - ix^2}  \frac{\zeta(s)\zeta(s-1)}{s(1-s)}x^s \,ds + O( x^\delta ).
   \label{eq:4.1.1}
\end{align}
By the same consideration for the rectangle 
\[
   \{ s = \sigma + it \mid 2 + \delta \leq \sigma \leq 3 + \delta,\ x^2 \leq t \leq T,\ x^2 <T \},
\]
we have
\begin{equation*}
   \int_{3 + \delta + ix^2}^{3 + \delta + i\infty} \frac{\zeta(s)\zeta(s-1)}{s(1-s)}x^s \,ds
   =  -\int_{2 + \delta + i x^2}^{3 + \delta + ix^2} \frac{\zeta(s)\zeta(s-1)}{s(1-s)}x^s \,ds  + O( x^\delta ).
\end{equation*}
Thus, we can rewrite the integral in \eqref{eq:3.1.12} along the curve $\mathcal{L}$ as follows:
\begin{equation}
\begin{aligned}
   &\left\{ \int_{3 + \delta - i\infty}^{3 + \delta - ix^2}
   + \int_{2 + \delta -ix^2}^{\frac{1}{2} + \eta -ix^2}
   + \int_{\frac{1}{2} + \eta -ix^2}^{\frac{1}{2} + \eta + ix^2}
   + \int_{\frac{1}{2} + \eta  + ix^2}^{ 2 + \delta + ix^2}
   + \int_{3 + \delta + ix^2}^{3 + \delta + i\infty} \right\}
\\[-4pt]
   &\hspace{2cm} \times \frac{\zeta(s)\zeta(s-1)}{s(1-s)}x^s \,ds,
\end{aligned}
   \label{eq:4.1.2}
\end{equation}
where $\eta = 1/\log\log x^2$. Moreover, the integral from $3 + \delta + ix^2$ to $3 + \delta + i\infty$ in \eqref{eq:4.1.2} is split as
\begin{equation}
   \left\{ \int_{2 + \delta + i\infty}^{2 + \delta + ix^2}
   + \int_{2 + \delta + i x^2}^{3 + \delta + ix^2} \right\}
   \frac{\zeta(s)\zeta(s-1)}{s(1-s)}x^s \,ds 
   = \int_{2 + \delta + i x^2}^{3 + \delta + ix^2}  \frac{\zeta(s)\zeta(s-1)}{s(1-s)}x^s \,ds + O( x^\delta ).
   \label{eq:4.1.3}
\end{equation}
The same estimate as \eqref{eq:4.1.3} holds for the integral from $3 + \delta - ix^2$ to $3 + \delta - i\infty$. From \eqref{eq:4.1.1} and \eqref{eq:4.1.2}, \eqref{eq:3.1.11} can be rewritten as follows:
\begin{align}
   \as 
   &= \frac{1}{2\pi i}
   \left\{
      \int_{2 + \delta -ix^2}^{\frac{1}{2} + \eta -ix^2}
      + \int_{\frac{1}{2} + \eta -ix^2}^{\frac{1}{2} + \eta +ix^2}
      +  \int_{\frac{1}{2} + \eta +ix^2}^{2 + \delta +ix^2}
   \right\}
   \frac{\zeta(s)\zeta(s-1)}{s(1-s)}x^s \,ds   \notag\\
   &\quad + O(x^{\delta + 3}).
   \label{eq:4.1.4}
\end{align}

Now we estimate each integral in \eqref{eq:4.1.4}. We use the following bound for $\log \zeta (s)$ under the Riemann Hypothesis.

\begin{fact}[~\cite{t}]
On the Riemann Hypothesis, we have
\begin{equation}
   \log \zeta(s) = O\left( \frac{(\log t)^{2-2\sigma}}{\log\log t} \right),
   \label{eq:4.1.5}
\end{equation}
uniformly for $1/2 < \sigma_0 \leq \sigma \leq \sigma_1 < 1$.
\end{fact}

Let $D$ be the rectangle
\[
   D = \left\{ s = \sigma + it \ \middle| \ \frac{1}{2} + \eta \leq \sigma \leq 2 + \delta,\ 0 < t_0 \leq t \leq x^2,\ t_0 \geq 1 \right\}.
\]First, we consider the case $\sigma = 1/2 + \eta$.Using Fact~4.1.1 and standard bounds for $\zeta(s)$, we have 
\[
   \zeta(s)\zeta(s-1)\ll |t|^{1 - \eta} \exp \left( \frac{(\log t)^{1-2\eta}}{\log\log t} \right)
   < |t| \exp \left( \frac{\log t}{\log\log t} \right)
   \ll |s| \exp \left( \frac{\log |s|}{\log\log |s|} \right),
\]
and hence
\begin{equation}
   \zeta(s)\zeta(s-1)\bigg/ \left( s\exp \left( \frac{\log |s|}{\log\log |s|} \right) \right) \ll 1
   \label{eq:4.1.6}
\end{equation}
for $\sigma = 1/2 + \eta$. Secondly, we consider the case $\sigma = 2 + \delta$. In the same way as in the derivation of \eqref{eq:4.1.6}, we obtain the same bound for $\sigma = 2 + \delta$. By the maximum modulus principle, we have 
\begin{equation}
   \zeta(s)\zeta(s-1) \ll |s|\exp\left( \frac{\log |s|}{\log \log |s|} \right)
   \label{eq:4.1.7}
\end{equation}
for all $s \in D$. In particular, for $s = \sigma + ix^2$ with $(1/2) + \eta \leq \sigma \leq 2 + \delta$, the left-hand side of \eqref{eq:4.1.7} satisfies
\begin{equation}
   \zeta(s)\zeta(s-1) \ll x^2 \exp\left( \frac{\log x}{\log\log x} \right).
   \label{eq:4.1.8}
\end{equation}
Using \eqref{eq:4.1.7}, the integral from $2+\delta -ix^2$ to $(1/2)+ \eta -ix^2$ in \eqref{eq:4.1.4} is bounded by
\begin{equation}
   \ll x^\delta \exp\left( \frac{\log x}{\log\log x} \right).
   \label{eq:4.1.9}
\end{equation}
Similarly, the same bound as \eqref{eq:4.1.9} holds for the integral from $(1/2)+\eta + ix^2$ to $2 + \delta + ix^2$ in \eqref{eq:4.1.4}. Now we estimate the integral from $1/2 + \eta - ix^2$ to $1/2 + \eta + ix^2$ in \eqref{eq:4.1.4}. For $s = (1/2)+\eta + it$ with $t_0 \leq |t| \leq x^2$, we have $2 - 2\sigma = -2\eta + 1$ by the definition of $\eta$. Using this and Fact~4.1.1, we obtain
\begin{equation}
   \zeta(s)\zeta(s-1)\ll |t| \exp \left( \frac{(\log |t|)^{-2\eta + 1}}{\log\log |t|} \right)
   \ll |t| \exp\left( \frac{\log x}{\log\log x} \right).
   \label{eq:4.1.10} 
\end{equation}
First, we consider the case $t_0 \leq t \leq x^2$. The integrand in \eqref{eq:4.1.4} from $(1/2) + \eta - ix^2$ to $(1/2) + \eta + ix^2$ is bounded by 
\begin{align}
   \zeta(s)\zeta(s-1)\frac{x^s}{s(1-s)}
   &\ll |t| \exp\left( \frac{\log x}{\log\log x} \right)
      \cdot \frac{x^{1/2 + \eta}}{|t|^2}   \notag\\
   &= \frac{x^{1/2 + \eta}}{|t|}\exp\left( \frac{\log x}{\log\log x} \right).
   \label{eq:4.1.11} 
\end{align}
Secondly, we consider the case $0<t\leq t_0$. Since we assume the Riemann Hypothesis, the Lindel\"of Hypothesis holds. Also, the complex number $s -1 = -1/2 + \eta + it$ tends to $-1/2 + it$ as $x \to \infty$, so $\zeta(s)$ satisfies 
\[
   \zeta(\sigma + it) \asymp |t|
\]
for $\sigma<0$. Hence, the left-hand side of \eqref{eq:4.1.11} is bounded by
\begin{equation}
   \ll \frac{x^{1/2+\eta}}{|t|^{1-\epsilon}} \ll x^{1/2 + \eta},   \label{eq:4.1.12}
\end{equation}
where $\epsilon$ is any positive real number. Therefore, the integral in question is estimated by 
\begin{align}
   \int_{\frac{1}{2} + \eta - ix^2}^{\frac{1}{2} + \eta + ix^2} \zeta(s)\zeta(s-1)\frac{x^s}{s(1-s)}\,ds
   &\ll x^{1/2 + \eta} + x^{1/2 + \eta} \exp\left( \frac{\log x}{\log\log x} \right) \int_{t_0}^{x^2} \frac{dt}{t}   \notag\\
   &= x^{1/2} \left\{ \exp\left( \frac{\log x}{\log\log x} \right) \right.   \notag\\
   &\quad \left.+ \exp\left( 2\frac{\log x}{\log\log x} \right) \exp(\log\log x^2) \right\}   \notag\\
   &\ll x^{1/2} \exp\left( \frac{\log x}{\log\log x} \right).
   \label{eq:4.1.13}
\end{align}
By \eqref{eq:4.1.13}, we obtain
\begin{align*}
   \as
   &\ll x^\delta \exp\left( \frac{\log x}{\log\log x} \right)
      + x^\delta \exp\left( \frac{\log x}{\log\log x} \right) \\
   &\quad +  x^{1/2} \exp\left( \frac{\log x}{\log\log x} \right) + O(x^{3+\delta}).
\end{align*}
Taking $\delta' = \max\{ 1/2 , \delta \}$, Theorem~1.2.1 follows. \qed

\subsection{Proof of Theorem 1.2.2}
We now prove Theorem 1.2.2. Let $\epsilon$ be any positive real number and $x \geq 1$. First, we consider the case $x > \exp (\exp \epsilon^{-1})$. Comparing the exponential factor in Theorem~1.2.1 with $x^\epsilon$ and using the condition on $x$, we have 
\[
   \exp\left( \frac{\log x}{\log \log x} \right) \ll x^\epsilon .
\]
By Theorem~1.2.1, we see that
\[
   \as \ll x^{\delta' + \epsilon},
\]
and we obtain the desired result for $x > \exp(\exp \epsilon^{-1})$. Secondly, we consider the case $1 \leq x \leq \exp (\exp \epsilon^{-1})$. We define the positive constant
\begin{equation}
   C = \max_{x \in [1, \exp(\exp \epsilon^{-1})]} \left| \as / x^{\delta' + \epsilon} \right|.
   \label{eq:4.2.1}
\end{equation}
Since $\as$ is absolutely continuous, the existence of the constant in \eqref{eq:4.2.1} is ensured. Hence, we have
\[
   \as \ll_\epsilon x^{\delta' + \epsilon}
\]
and obtain the desired result for $1 \leq x \leq \exp(\exp \epsilon^{-1})$. Therefore, Theorem 1.2.2 is proved. \qed

\subsection{Why the RH-equivalent condition disappers}
We imitate the proof of the Riemann Hypothesis under the bound \eqref{eq:1.1.4} in Theorem 1.1.1. As in the proof of the Mellin transform for $\as$ in \eqref{eq:3.1.1}, the following Mellin transform for $\an$ holds:
\begin{equation}
   \int_{1}^\infty \an x^{-s-1}\,dx
   = \frac{3}{\pi^2}\frac{1}{s-2} + \frac{1}{s(1-s)}\frac{\zeta(s-1)}{\zeta(s)}
   \label{eq:4.2.2}
\end{equation}
(Lemma~5.2 in~\cite{Kac and Wie}). Let $\Delta$ be any compact subset of the half-plane $\text{Re}(s) > 1/2$. Also, let $d = \inf_{s \in \Delta} \text{Re} (s)$ with $\epsilon < d$. For $s \in \Delta$, we see that $\text{Re} (s) \geq 1/2 + d$ and
\[
   \int_{1}^\infty \an x^{-s-1}\,dx
   \ll_\epsilon \int_{1}^\infty x^{1/2 + \epsilon} \cdot x^{-3/2 - d}\,dx
   = \int_{1}^\infty x^{-1 + \epsilon - d}\,dx
   = \frac{1}{d - \epsilon}.
\]
Hence, the integral on the left-hand side is holomorphic in the region $\text{Re}(s) > 1/2$ and $\zeta(s) \neq 0$ there. Using the functional equation for $\zeta(s)$, the same conclusion holds in the region $\text{Re} (s) < 1/2$. Therefore, the Riemann Hypothesis holds and so Theorem~1.1.1 holds. 

The above proof works because there is the factor $\zeta(s)^{-1}$ on the right-hand side of \eqref{eq:4.2.2}.The factor $\zeta(s-1)/\zeta(s)$ appears from the Dirichlet series expansion
\begin{equation}
   \sum_{n=1}^\infty \frac{\varphi(n)}{n^s} = \frac{\zeta(s-1)}{\zeta(s)}
   \label{eq:4.2.3}
\end{equation}for $\sigma > 2$.It seems that Theorem~1.1.1 can be obtained whenever the factor $\zeta(s)$ appears in the denominator of the Dirichlet series for $b(n)$ in Theorem 2.1.1.In~\cite{Kac and Wie}, the case where $b(n) = \varphi(n)$ in Theorem~2.1.1 is considered, and Theorem~1.1.1 is proved as an equivalent condition to the Riemann Hypothesis. However, we cannot obtain such an equivalence to the Riemann Hypothesis when we consider the case $b(n) = \sigma_1(n)$ in Theorem~2.1.1. It seems that the same phenomenon occurs for other arithmetic functions of the type of the divisor function $\sigma_1 (n)$ in Theorem~2.1.1.

%\textbf{Acknowledgments.} The author appreciates valuable comments from the anonymous referee. In particular, the author referred to comments from the referee about the remark of Section~2.  

\end{document}